\documentclass[11pt]{article}

\usepackage{amssymb}

\newcommand{\newc}{\newcommand}

\newc{\eqnoset}{\setcounter{equation}{0}}

\newcommand{\mref}[1]{(\ref{#1})}
\newcommand{\reflemm}[1]{Lemma~\ref{#1}}
\newcommand{\refrem}[1]{Remark~\ref{#1}}
\newcommand{\reftheo}[1]{Theorem~\ref{#1}}
\newcommand{\refdef}[1]{Definition~\ref{#1}}

\newcommand{\refsec}[1]{Section~\ref{#1}}

\newcommand{\beq}{\begin{equation}}
\newcommand{\eeq}{\end{equation}}
\newcommand{\beqno}[1]{\begin{equation}\label{#1}}

\newcommand{\barr}{\begin{array}}
\newcommand{\earr}{\end{array}}

\newc{\bearr}{\begin{eqnarray*}}
\newc{\eearr}{\end{eqnarray*}}

\newc{\bearrno}[1]{\begin{eqnarray}\label{#1}}
\newc{\eearrno}{\end{eqnarray}}

\newc{\non}{\nonumber}
\newc{\nol}{\nonumber\nl}

\newcommand{\bdes}{\begin{description}}
\newcommand{\edes}{\end{description}}
\newc{\benu}{\begin{enumerate}}
\newc{\eenu}{\end{enumerate}}
\newc{\btab}{\begin{tabular}}
\newc{\etab}{\end{tabular}}

\newtheorem{theorem}{Theorem}[section]
\newtheorem{defi}[theorem]{Definition}
\newtheorem{lemma}[theorem]{Lemma}
\newtheorem{rem}[theorem]{Remark}
\newtheorem{exam}[theorem]{Example}
\newtheorem{propo}[theorem]{Proposition}
\newtheorem{corol}[theorem]{Corollary}

\newcommand{\btheo}[1]{\begin{theorem}\label{#1}}
\newc{\brem}[1]{\begin{rem}\label{#1}\em}
\newc{\bexam}[1]{\begin{exam}\label{#1}\em}
\newc{\bdefi}[1]{\begin{defi}\label{#1}}
\newcommand{\blemm}[1]{\begin{lemma}\label{#1}}
\newcommand{\bprop}[1]{\begin{propo}\label{#1}}
\newcommand{\bcoro}[1]{\begin{corol}\label{#1}}
\newcommand{\etheo}{\end{theorem}}
\newcommand{\elemm}{\end{lemma}}
\newcommand{\eprop}{\end{propo}}
\newcommand{\ecoro}{\end{corol}}
\newc{\erem}{\end{rem}}
\newc{\eexam}{\end{exam}}
\newc{\edefi}{\end{defi}}

\newc{\rmk}[1]{{\bf REMARK #1: }}
\newc{\DN}[1]{{\bf DEFINITION #1: }}

\newcommand{\bproof}{{\bf Proof:~~}}
\newc{\eproof}{{\vrule height8pt width5pt depth0pt}\vspace{3mm}}

\newc{\bfrac}[2]{\dspl{\frac{#1}{#2}}}

\newc{\nid}{\noindent}

\newcommand{\dspl}{\displaystyle}
\newc{\grad}{\nabla}
\newc{\Div}{\mbox{div}}
\newc{\pdt}[1]{\dspl{\frac{\partial{#1}}{\partial t}}}
\newc{\pdn}[1]{\dspl{\frac{\partial{#1}}{\partial \nu}}}
\newc{\pdNi}[1]{\dspl{\frac{\partial{#1}}{\partial \mathcal{N}_i}}}
\newc{\pD}[2]{\dspl{\frac{\partial{#1}}{\partial #2}}}
\newc{\dt}{\dspl{\frac{d}{dt}}}
\newc{\bdry}[1]{\mbox{$\partial #1$}}
\newc{\sgn}{\mbox{sign}}

\newc{\Hess}[1]{\frac{\partial^2 #1}{\pdh z_i \pdh z_j}}
\newc{\hess}[1]{\partial^2 #1/\pdh z_i \pdh z_j}

\newc{\ag}{\alpha}
\newc{\bg}{\beta}
\newc{\cg}{\gamma}\newc{\Cg}{\Gamma}
\newc{\dg}{\delta}\newc{\Dg}{\Delta}
\newc{\eg}{\varepsilon}
\newc{\zg}{\zeta}
\newc{\thg}{\theta}
\newc{\llg}{\lambda}\newc{\LLg}{\Lambda}
\newc{\kg}{\kappa}
\newc{\rg}{\rho}
\newc{\sg}{\sigma}\newc{\Sg}{\Sigma}
\newc{\tg}{\tau}
\newc{\fg}{\phi}\newc{\Fg}{\Phi}
\newc{\vfg}{\varphi}
\newc{\og}{\omega}\newc{\Og}{\Omega}
\newc{\pdh}{\partial}

\newc{\ccG}{{\cal G}}

\newc{\ii}[1]{\int_{#1}}
\newc{\iidx}[2]{{\dspl\int_{#1}~#2~dx}}
\newc{\bii}[1]{{\dspl \ii{#1} }}
\newc{\biii}[2]{{\dspl \iii{#1}{#2} }}
\newc{\su}[2]{\sum_{#1}^{#2}}
\newc{\bsu}[2]{{\dspl \su{#1}{#2} }}

\newc{\biiom}[1]{{\dspl\int_{\bdrom}~ #1 ~d\sg}}
\newc{\io}[1]{{\dspl\int_{\Og}~ #1 ~dx}}
\newc{\bio}[1]{{\dspl\int_{\bdrom}~ #1 ~d\sg}}
\newc{\bsir}{\bsu{i=1}{r}}
\newc{\bsim}{\bsu{i=1}{m}}

\newc{\iibr}[2]{\iidx{\bprw{#1}}{#2}}
\newc{\Intbr}[1]{\iibr{R}{#1}}
\newc{\intbr}[1]{\iibr{\rg}{#1}}
\newc{\intt}[3]{\int_{#1}^{#2}\int_\Og~#3~dxdt}

\newc{\itQ}[2]{\dspl{\int\hspace{-2.5mm}\int_{#1}~#2~dz}}
\newc{\mitQ}[2]{\dspl{\rule[1mm]{4mm}{.3mm}\hspace{-5.3mm}\int\hspace{-2.5mm}\int_{#1}~#2~dz}}
\newc{\mitQQ}[3]{\dspl{\rule[1mm]{4mm}{.3mm}\hspace{-5.3mm}\int\hspace{-2.5mm}\int_{#1}~#2~#3}}

\newc{\mitx}[2]{\dspl{\rule[1mm]{3mm}{.3mm}\hspace{-4mm}\int_{#1}~#2~dx}}
\newc{\mitmu}[2]{\dspl{\rule[1mm]{3mm}{.3mm}\hspace{-4mm}\int_{#1}~#2~d\mu}}
\newc{\iidmu}[2]{{\dspl\int_{#1}~#2~d\mu}}

\newc{\iidm}[3]{{\dspl\int_{#1}~#2~d #3}}

\newc{\itQmu}[2]{\dspl{\int\hspace{-2.5mm}\int_{#1}~#2~d\mu}}
\newc{\mitQmu}[2]{\dspl{\rule[1mm]{4mm}{.3mm}\hspace{-5.3mm}\int\hspace{-2.5mm}\int_{#1}~#2~d\mu}}

\newc{\mitQq}[2]{\dspl{\rule[1mm]{4mm}{.3mm}\hspace{-5.3mm}\int\hspace{-2.5mm}\int_{#1}~#2~d\bar{z}}}
\newc{\itQq}[2]{\dspl{\int\hspace{-2.5mm}\int_{#1}~#2~d\bar{z}}}

\newc{\pder}[2]{\dspl{\frac{\partial #1}{\partial #2}}}

\newc{\bdrom}{\bdry{\Og}}

\newc{\bilhom}{\mbox{Bil}(\mbox{Hom}(\RR^{nm},\RR^{nm}))}
\newc{\VV}[1]{{V(Q_{#1})}}

\newc{\ccA}{{\mathcal A}}
\newc{\ccB}{{\mathcal B}}
\newc{\ccC}{{\mathcal C}}
\newc{\ccD}{{\mathcal D}}
\newc{\ccE}{{\mathcal E}}
\newc{\ccH}{\mathcal{H}}
\newc{\ccF}{\mathcal{F}}
\newc{\ccI}{{\mathcal I}}
\newc{\ccJ}{{\mathcal J}}
\newc{\ccK}{{\mathcal K}}
\newc{\ccP}{{\mathcal P}}
\newc{\ccQ}{{\mathcal Q}}
\newc{\ccR}{{\mathcal R}}
\newc{\ccS}{{\mathcal S}}
\newc{\ccT}{{\mathcal T}}
\newc{\ccX}{{\mathcal X}}
\newc{\ccY}{{\mathcal Y}}
\newc{\ccZ}{{\mathcal Z}}

\newc{\bb}[1]{{\mathbf #1}}

\newc{\myprod}[1]{\langle #1 \rangle}
\newc{\mypar}[1]{\left( #1 \right)}

\newc{\BLLg}{\mathbf{\LLg}}

\newc{\mA}{\mathbf{A}}
\newc{\mB}{\mathbf{B}}
\newc{\mC}{\mathbf{C}}
\newc{\mD}{\mathbf{D}}
\newc{\mE}{\mathbf{E}}
\newc{\mF}{\mathbf{F}}
\newc{\mJ}{\mathbf{J}}
\newc{\mG}{\mathbf{G}}
\newc{\mP}{\mathbf{P}}
\newc{\mR}{\mathbf{R}}
\newc{\mQ}{\mathbf{Q}}
\newc{\mX}{\mathbf{X}}
\newc{\muu}{\mathbf{u}}
\newc{\mvv}{\mathbf{v}}

\newc{\mllg}{\mathbb{\lambda}}
\newc{\mLLg}{\mathbf{\LLg}}

\newc{\lspn}[2]{\mbox{$\| #1\|_{\Lsp{#2}}$}}
\newc{\Lpn}[2]{\mbox{$\| #1\|_{#2}$}}
\newc{\Hn}[1]{\mbox{$\| #1\|_{H^1(\Og)}$}}

\newc{\mynorm}[2]{\| #1\|_{#2}}

\newcommand{\RR}{{\rm I\kern -1.6pt{\rm R}}}

\newc{\itQQ}[2]{\dspl{\int_{#1}#2\,dz}}
\newc{\mmitQQ}[2]{\dspl{\rule[1mm]{4mm}{.3mm}\hspace{-4.3mm}\int_{#1}~#2~dz}}
\newc{\MmitQQ}[2]{\dspl{\rule[1mm]{4mm}{.3mm}\hspace{-4.3mm}\int_{#1}~#2~d\mu}}

\newc{\MUmitQQ}[3]{\dspl{\rule[1mm]{4mm}{.3mm}\hspace{-4.3mm}\int_{#1}~#2~d#3}}
\newc{\MUitQQ}[3]{\dspl{\int_{#1}~#2~d#3}}

\newc{\mccP}{\mathbb{P}}
\newc{\mccK}{\mathbb{K}}

\newc{\DKTmU}{\mccK(U)}
\newc{\DKTmUold}{(K_U(U)^{-1})^T}

\newc{\myPi}{\mathbf{W}}
\newc{\myIbar}{\bar{\ccI}_1}
\newc{\myIhat}{\hat{\ccI}_1}
\newc{\myIbreve}{\breve{\ccI}_0}

\newc{\mmk}{\mathbf{k}}

\newcommand{\ma}{\mathbf{a}}
\newcommand{\mg}{\mathbf{g}}

\newc{\mfu}{\mathbf{f_u}}
\newc{\mh}{\mathbf{h}}
\newc{\mb}{\mathbf{b}}

\newcommand{\barrl}[2]{\barr{ll}\lefteqn{#1}\hspace{#2}&\\}

\newc{\mN}{\mathbf{N}}
\newc{\mI}{\mathbf{I}}
\newc{\mH}{\mathbf{H}}

\newc{\mk}{\mathbf{k}}
\newc{\mr}{\mathbf{r}}

\newc{\DIAGM}[2]{\left[\barr{ccc}#1&0\ldots&0\\
	\vdots&\ddots&\vdots\\0&\ldots0&#2\earr \right]}
\newc{\DiagM}[2]{\mbox{diag}\left[#1
	\cdots #2 \right]}
\newc{\vVEC}[2]{\left[\barr{c}#1\\
	\vdots\\#2\earr \right]}
\newc{\hVEC}[2]{\left[#1
	\cdots #2 \right]}

\newc{\mq}{\mathbf{q}}

\newc{\msys}[1]{\left\{\barr{l}#1\earr
	\right.}
\newc{\msysa}[1]{\left\{\barr{ll}#1\earr
	\right.}

\newc{\bbM}{\mathbb{M}}
\newc{\mat}[1]{\left[\barr{cc}#1\earr\right]}

\begin{document}

\vspace*{-.8in}
\begin{center} {\LARGE\em On the Global Existence of a Class of Strongly Coupled Parabolic Systems.}

 \end{center}

\vspace{.1in}

\begin{center}

{\sc Dung Le}{\footnote {Department of Mathematics, University of
Texas at San
Antonio, One UTSA Circle, San Antonio, TX 78249. {\tt Email: Dung.Le@utsa.edu}
}}

\end{center}

\begin{abstract}
We establish the global existence of a class of strongly coupled parabolic  systems. The necessary apriori estimates will be obtained via our new approach to the regularity theory of parabolic scalar equations with integrable data and new $W^{1,p}$ estimates of their solutions. The key assumption here is that the $L^p$ norms of solutions are uniformly bounded for some sufficiently large $p\in (1,\infty)$, an assumption can be easily affirmed for systems with polynomial growth data. This replaces the usual condition that the solutions are uniformly bounded which is very hard to be verified because maximum principles for systems are generally unavailable. \end{abstract}

\vspace{.2in}

\section{Introduction} \label{intro}\eqnoset

In this paper, let $\Og$ be a bounded domain in $\RR^N$, $N\ge2$, with smooth boundary $\partial \Og$ and $T$ be a positive number.  We study the solvability of the strongly coupled parabolic system
\beqno{fullsys0a}\msysa{W_t=-\Div(\mA DW)+\mB DW+\mG W+F&\mbox{in $\Og\times(0,T)$,}\\ \mbox{Homogeneous Dirichlet or Neumann boundary conditions}&\mbox{on $\partial\Og$,}\\W=W_0&\mbox{on $\Og\times(0,T)$.}}\eeq
Here,  $W=[u_i]_{i=1}^m$, a vector in $\RR^m$.  $\mA,\mB,\mG$ are square matrices of size $m\times m$ and $F=[f_i]_{i=1}^m$ for some $m\ge2$. The entries of these matrices are functions in $W$. $W_0$ is a given vector valued function in $W^{1,N_0}(\Og)$ with some $N_0>N$.

In particular, the entries of $\mB$ are vectors in $\RR^N$ and the entries of $\mB DW$ should be understood as $\sum \myprod{\mB_{ij},Du_i}$.

The existence problem of the system \mref{fullsys0a} was investigated by Amann. He shows that if the parameters of the regular parabolic \mref{fullsys0a} (i.e. $\mA$ is normal elliptic (see \cite{Am2})) are bounded and  \beqno{est1a} \limsup_{t\to T}\|W\|_{W^{1,2p}(\Og)} <\infty \mbox{ for some $p>N/2$}\eeq then there is a unique strong solution $W\in W^{2,2}_{loc}(\Og\times (0,T))$ of \mref{fullsys0a}, with $DW$ bounded and satisfies the system a.e. in $\Og\times (0,T)$.

For nonlinear strongly coupled systems like \mref{fullsys0a} one would start by proving the boundedness of solutions because Amann's theory worked with bounded $\mA$. For nonlinear strongly coupled systems like \mref{fullsys0a} this problem would be a very hard one already (see counterexamples in \cite{JS, LN}). The next obstacle, and even harder, is the estimate of higher order norms like \mref{est1a}. This problem is closely related to the regularity of parabolic systems in \cite{bj,Cam,DM,GiaS}.

In this paper, we will suppose that certain integrabilities of solutions are available, an assumption can be verified in important cases (e.g. see \cite{SKT,yag}), in order to deduce \mref{est1a}. To this end, we will consider each equation in the system \mref{fullsys0a} and take a close look at the regularity of nonlinear scalar equations with integrable coefficients. We will extend the idea in \cite{dleNA} to the case when  only certain integrabilities of the data are available. Note that even when we establish the boundedness of solutions themselves {\em we cannot suppose that these data are bounded because they depend on the components of the other equations and whose boundedness are not known yet}. Importantly, $W^{1,2p}$ estimates like \mref{est1a} are essentially needed due to the strong couplings of \mref{fullsys0a}. Therefore, the regularity theory for scalar equations studied in \refsec{scalareqn}, besides its generality and its own interests, needs closer looks at this matter.

We collect technical known facts and their improvements in \refsec{techlem} for their uses later. One  important ingredient is a compactness result of a sequence of solutions to a family of \mref{fullsys0a} in $L^2((0,T), L^p(\Og))$. This is a consequence of the Simon-Aubin-Lions compact lemma. We present the details and its  consequence: weighted interpolation inequality \mref{interpol}, an important tool in the proofs discussed in \refsec{scalareqn}.

In stating the results and presenting their proofs, we introduce the following definition
\bdefi{Ogdef} We say that a function $f:\Og\times(0,T)\to \RR$ is of class $\bbM(\Og,T)$ if there is some $p_0>N/2$ such that (the number $p_0$ is not fixed) \beqno{Mnorm}\sup_{(0,T)}\|f\|_{L^{p_0}(\Og)}<\infty.\eeq\edefi
This definition also applies to matrix valued functions. We say that a matrix is in $\bbM(\Og,T)$ if its entries are.

In \refsec{scalareqn}, we will consider the scalar equations in $Q=\Og\times(0,T)$ written in the forms
\beqno{diagSKTa} v_t=\Div(ADv)+\Div(B_1)+B_2Dv+G_1v+G_ ,\eeq
\beqno{eq1a}v_t=\Div(\ma Dv)+\hat{\mb}Dv+\Div(\mb)+\mg,\eeq
with homogeneous Neumann or Dirichlet conditions on $\partial \Og\times(0,T)$.
Here $B_1,B_2,\hat{\mb},\mb$ are row vectors in $\RR^N$ and $A,\ma,G_1,G_2,\mg$ are scalar functions in $\bbM(\Og,T)$.

We will show that the weak solutions of these equations are bounded and H\"older continuous under the assumption that the parameters of the equations are in $\bbM(\Og,T)$. What is important here is that their norms (and H\"older exponents) are uniform with respect to the norms of their parameters in $\bbM(\Og,T)$. Furthermore, combining with the results in \cite{dleANS},  we also present the $W^{1,2p}(\Og)$ estimates of these equations, a crucial fact used in the next section where we consider strongly coupled systems and the ultimate goal is to obtain \mref{est1a}.

In \refsec{sys}, we will apply the theory of \refsec{scalareqn} to each equation of \mref{fullsys0a} in order to derive \mref{est1a}. Examples are provided to see that the condition that the parameters of \mref{fullsys0a} are in $\bbM(\Og,T)$ can be verified for a class of important systems in applications.

We consider first triangular systems \mref{fullsys0a}
where  $\mA,\mB$ are upper triangular matrices and $\mG=[g_{ij}]$ a full matrix. We will prove (by induction) that
\btheo{trisysthma} Assume that $\mA=[\ag_{ij}]$, $\mB=[\bg_{ij}]$ are upper triangular matrices ($\ag_{ij}=\bg_{ij}=0$ if $i<j$) and $\mG=[g_{ij}]$ is a full matrix and $F=[f_i]$. 
In particular, we will assume that $\ag_{ij}$ depends only on $u_j$ for $j\ge i$. Furthermore, assume the following integrability conditions
\beqno{trisys1} \ag_{ij}, f_{i}\in \bbM(\Og,T).\eeq

Suppose that for some $p>N/2$ and any $i<j$ \beqno{triFsys}|(\bg_{ij})_{u_i}Du_i|\in L^\frac{2p+2}{2}(Q) \mbox{ and }|(g_{ij})_{u_i}Du_i|,\;|(f_j)_{u_i}Du_i|\in L^\frac{2p+2}{3}(Q).\eeq

And for some $q$ such that $p>q>N/2$ and any $i<j$
\beqno{trisysbga} \iidx{\Og}{\ag_{i,j}^\frac{2qp}{p-q}},\; \iidx{\Og}{\bg_{i,j}^\frac{2qp}{2p-q}}< \infty. \eeq

Then \mref{fullsys0a} has a unique strong solution on $(0,T)$ for any $N$.
\etheo

We then move on to the case of \mref{fullsys0a} being a full system. A simple observation reveals that we need to establish a weaker version of \mref{est1a} in general. That is, we need only to establish a uniform estimate for the $W^{1,2p}$ norms of some components of $W$. To illustrate this, let us state the result when $m=2$. 

Consider \mref{fullsys0a} with $\mA=\mat{\ag&\bg\\\cg&\dg}, \mB=\mat{a&b\\c&d}$, $\mG=\mat{g_1&0\\0&g_2}$ are full $2\times2$ matrices and $\mF=[f_i]$. Assume that the entries of theses matrices are smooth functions in $u,v$. Under some appropriate integrability assumptions  of these parameters, we will show that

\btheo{fullthma} Consider \mref{fullsys0a} when $m=2$ with   $W=[u,v]^T$. Suppose that for  some $p>N/2$  we can establish a bound for $\sup_{(0,T)}\|u\|_{W^{1,2p}(\Og)}$ with $2p>N$ (compare with \mref{est1a} which requires also a bound for $\sup_{(0,T)}\|v\|_{W^{1,2p}(\Og)}$) and 
\beqno{fulbbM}\dg,cDu,d^2,g_2\in\bbM(\Og,T),\;\eeq
\beqno{triFa2}|d_u Du|\in  L^\frac{2p+2}{2}(Q)\mbox{ and }|(g_2)_uDu|,\;|(f_2)_uDu|\in L^\frac{2p+2}{3}(Q).\eeq

 Suppose also that for some $q$ such that $p>q>N/2$ and
$$ \iidx{\Og}{\cg^\frac{2qp}{p-q}}  <\infty.$$

Then \mref{fullsys0a} has a unique strong solution on $(0,T)$. \etheo

The assumption on $\sup_t\|u\|_{W^{1,2p}(\Og)}$ of this result is clearly weaker than \mref{est1a} of Amann's which need the estimates for both components of $W$. Of course, simple uses of Young's inequality and the assumption that $\sup_{(0,T)}\|u\|_{W^{1,2p}(\Og)}$ is bounded would show that the condition on $cDu$ in \mref{fulbbM}  and \mref{triFa2} (which can be dropped, see \refrem{Gap} at the end of the paper) can be implied from $$\sup_t|c|\in L^2(\Og), d_u\in L^\frac{2p}{p-1}(Q) \mbox{ and } g_u, (f_2)_u\in L^\frac{3p}{2p-1}(Q).$$

\section{Some technical lemmas}\label{techlem}\eqnoset

Similar to \cite[Lemma 3.3]{dletrans}, we will establish the following compactness result which will lead us to a key interpolation inequality in \reflemm{Interpolation}.

\blemm{complem} Suppose that $\llg_0$ is a positive constant and  and $\Fg_0,\Fg\in L^1(Q)$ such that $\llg_0\le \Fg_0\le\Fg$. Let $\ccF$ be a collection of function $v$ satisfying $\itQ{Q}{\Fg(|v|^2+|Dv|^2)}dt <M$ for some constant $M$ and
\beqno{vccF}\left|\itQ{Q}{v_t\psi}\right|\le C\itQ{Q}{(\Fg_0|Dv||D\psi|+\Fg|v||\psi|)} \quad \forall\psi\in C^1(Q).\eeq

If  $p$ is such that   $W^{1,2}(\Og)$  is compactly embedded in $\ccB := L^p(\Og)$  then $\ccF$ is compactly embedded in $L^2((0,T),\ccB)$.
\elemm

\bproof  We first show that for $l>(N+2)/2$ \beqno{wcont}\|v(\cdot,t+h)-v(\cdot,t)\|_{W^{-l,2}(\Og)}\le O(h).\eeq

For any $Q_{s,t}=\Og\times(s,t)$ let $\psi(x,t)=\eta(t)\fg(x)$ with $\fg\in C^1(\Og)$ and $\eta\equiv1$ in $(s,r)$ and $\eta\equiv0$ outside $(s-\eg,r+\eg)$. From \mref{vccF} and the assumption on $\Fg_0,\Fg$, we get
$$\barr{lll}\left|\itQ{Q_{s,t}}{v_t\psi}\right|&\le& C\itQ{Q_{s,t}}{|\Fg\eta(|Dv||D\fg|+|v||\fg|)|}\\&\le&
\|\Fg\|_{L^1(Q_{s,t})}^\frac12
\left(\itQ{Q}{\Fg(|D v|^2+|v|^2)}\right)^\frac12\|\fg\|_{C^1(\Og)}\le O(|h|)\|\fg\|_{C^1(\Og)}.\earr$$

We used the fact that $\Fg\in L^1(Q)$,  by H\"older's inequality $\|\Fg\|_{L^1(Q_{s,t})}\le O(|h|)$ with $h=s-t$. Because $l>(N+2)/2$, $\|\fg\|_{C^1(\Og)}\le C\|\fg\|_{W^{l,2}(\Og)}$, this implies
$$\left|\itQ{Q_{s,t}}{v_t\psi}\right|\le O(|h|)\|\fg\|_{W^{l,2}(\Og)} \quad \forall\fg\in W^{l,2}(\Og).$$  Since $v\psi=0$ at $s,t$, we have $\itQ{Q_{s,t}}{v\psi_t}=-\itQ{Q_{s,t}}{v_t\psi}$. Letting $\eg\to0$, we also have
$$\left|\iidx{\Og}{[v(\cdot,s)-v(\cdot,t)]\fg}\right|=\left|\itQ{Q_{s,t}}{v\psi_t}\right|\le O(|h|)\|\fg\|_{W^{l,2}(\Og)} \quad \forall\fg\in W^{l,2}(\Og).$$

Hence \mref{wcont} follows (compare with \cite[Lemma 3.2]{dletrans}). We now interpolate $\ccB\cap W^{1,2}(\Og)$ between $W^{1,2}(\Og)$ and $W^{-l,2}(\Og)$ to get for any $\mu>0$ and $v\in \ccB\cap W^{1,2}(\Og)$ that
$$\|v(\cdot,t+h)-v(\cdot,t)\|_\ccB^2 \le \mu \|v\|_{W^{1,2}(\Og)}^2+C(\mu)\|v(\cdot,t+h)-v(\cdot,t)\|_{W^{-l,2}(\Og)}^2$$

Since $\dspl{\int_0^T}\|v\|_{W^{1,2}(\Og)}^2dt\le C(M)$ if $v\in \ccF$, for any given $\eg>0$, we can choose $\mu$ small first and then $h$ small such that by \mref{wcont}
$$\int_{0}^{T-h}\|v(\cdot,t+h)-v(\cdot,t)\|_\ccB^2dt<\eg.$$
Thus, we just prove that
$$\int_{0}^{T-h}\|v(\cdot,t+h)-v(\cdot,t)\|_\ccB^2dt\le O(h).$$

For any $t_1,t_2\in (0,T)$ and $v\in \ccF$, the fact that the collection  $\int_{t_1}^{t_2}vdt$ is a pre-compact set of $\ccB$ is clear because the set $$\{ \frac{1}{|t_2-t_1|}\int_{t_1}^{t_2}vdt\,:\,  \int_{t_1}^{t_2}\iidx{\Og}{\Fg(|v|^2+|Dv|^2)}dt <M)\}$$ belongs to the closure of the convex hull in $\ccB$, a bounded set in $W^{1,2}(\Og)$ as $\Fg\ge \llg_0>0$, which is compact in $\ccB$.

By a result of Simon \cite{JS} as in \cite[Lemma 3.3]{dletrans}, we proved the compactness of $\ccF$ in $L^2(0,T,\ccB)$. \eproof

\brem{comprem} In some cases we need to modify the equation defining $\ccF$ \mref{vccF} can be modified by adding an extra term $\hat{\Fg}|v||D\psi|$ in \mref{vccF} as follows
\beqno{vccF1}\left|\itQ{Q}{v_t\psi}\right|\le C\itQ{Q}{[\Fg_0|Dv||D\psi|+\Fg|v||\psi|+\hat{\Fg}|v||D\psi|]} \quad \forall\psi\in C^1(Q).\eeq
If $\Fg, \frac{\hat{\Fg}^2}{\Fg}\in L^1(Q)$ then we still have $\ccF$ is compactly emmbedded in $L^2((0,T),\ccB)$.

The proof can go on as before as long as we can establish \mref{wcont}. To this end, we just need to apply H\"older's inequality to the extra term $\hat{\Fg}|v|$
$$\barr{lll}\left|\itQ{Q_{s,t}}{v_t\psi}\right|&\le& C\itQ{Q_{s,t}}{[|\Fg\eta(|Dv||D\fg|+|v||\fg|)|+\hat{\Fg}|v||D\psi|]}\\&\le&
[\|\Fg\|_{L^1(Q_{s,t})}^\frac12+\|\frac{\hat{\Fg}^2}{\Fg}\|_{L^1(Q_{s,t})}^\frac12]
\left(\itQ{Q}{\Fg(|D v|^2+|v|^2)}\right)^\frac12\|\fg\|_{C^1(\Og)}.\earr$$

The factors $\|\Fg\|_{L^1(Q_{s,t})}^\frac12, \|\frac{\hat{\Fg}^2}{\Fg}\|_{L^1(Q_{s,t})}^\frac12\le O(h)$ because we are assuming  $\Fg, \frac{\hat{\Fg}^2}{\Fg}\in L^1(Q)$. We obtain \mref{wcont} again. The proof then continues.

\erem

We now establish the following interpolation inequality which will play an essential role in many proofs. One should note that \mref{interpol} below is not true in general without \mref{keyccF} unless $\Fg,\Fg_0$ are positive constants.

\blemm{Interpolation} Let $\llg_0,\Fg_0,\Fg$ be as \reflemm{complem}. Assume that $v$ satisfies
\beqno{keyccF}\left|\itQ{Q}{v_t\psi}\right|\le C\itQ{Q}{(\Fg_0|Dv||D\psi|+\Fg|v||\psi|)} \quad \forall\psi\in C^1(Q).\eeq
Suppose also that for some $p_0>N/2$ we have \beqno{keyFg}\Fg\in L^1(Q)\mbox{ and }\sup_{t\in (0,T)}\iidx{\Og}{\Fg^{p_0}}<\infty.\eeq

Then for any $\eg>0$ and $\bg\in(0,2]$ there is a finite constant $C(\eg,\bg)$ such that 
\beqno{interpol}\itQ{Q}{\Fg |v|^2}\le \eg \itQ{Q}{\Fg_0 |Dv|^2}+C(\eg,\bg)\|\Fg\|_{L^1(Q)}\left(\itQ{Q}{ |v|^\bg}\right)^{2/\bg}.\eeq
\elemm

\bproof Suppose that this is not true then there are $\eg_0>0$ and a sequence $\{v_n\}$ such that
$$\itQ{Q}{\Fg |v_n|^2}\ge \eg_0 \itQ{Q}{\Fg_0|Dv_n|^2}+n\|\Fg\|_{L^1(Q)}\left(\itQ{Q}{ |v_n|^\bg}\right)^{2/\bg}.$$

By scaling, we can assume that $\itQ{Q}{\Fg|v_n|^2}=1$ so that
$$1\ge \eg_0 \itQ{Q}{\Fg_0|Dv_n|^2}+n\|\Fg\|_{L^1(Q)}\left(\itQ{Q}{ |v_n|^\bg}\right)^{2/\bg}.$$
Let $\Fg_0=\llg_0$ this implies $\{v_n\}\subset \ccF$ for $M=1+\frac{1}{\eg_0}$ and $$\lim_{n\to\infty}\itQ{Q}{ |v_n|^\bg}=0.$$

Let $p>2$ be such that $\frac{p}{p-2}=p_0>N/2$. We have $p\in(2,\frac{2N}{N-2})$ so that $W^{1,2}(\Og)$ is compactly embedded in $\ccB=L^p(\Og)$. By the compactness of \reflemm{complem}, we can assume that $v_n\to v_*$ in $L^2((0,T),L^p(\Og))$. The above limit implies that  $\itQ{Q}{ |v_*|^\bg}=0$. The uniqueness of limits yields
$v_*=0$.  So, $v_n\to 0$ in $L^2((0,T),L^p(\Og))$. Therefore, $\|v_n\|_{L^2((0,T),L^p(\Og))}\to0$.

But for  $C_1=\sup_t\left(\iidx{\Og}{\Fg^\frac{p}{p-2}}\right)^{1-\frac{2}{p}}<\infty$ (by the assumption \mref{keyFg})
$$1=\itQ{Q}{\Fg |v_n|^2}\le \int_0^T\left(\iidx{\Og}{\Fg^\frac{p}{p-2}}\right)^{1-\frac{2}{p}}\left(\iidx{\Og}{|v_n|^{p}}\right)^\frac{2}{p} dt\le C_1\|v_n\|_{L^2((0,T),\ccB)}^2.$$

Hence, $\|v_n\|_{L^2((0,T),\ccB)}^2\ge 1/C_1$ for all $n$. This is a contradiction to $\|v_n\|_{L^2((0,T),L^p(\Og))}\to0$. The dependence of $C$ in \mref{interpol} on $\|\Fg\|_{L^1(Q)}$ can be seen by taking $v\equiv1$ on $Q$. The lemma is proved. \eproof

\brem{intrem}  Assume that $\Fg, \frac{\hat{\Fg}^2}{\Fg}\in L^1(Q)$  and for some $p_0>N/2$ \beqno{keyFg1} \sup_{t\in (0,T)}\iidx{\Og}{\Fg^{p_0}} \mbox{ and } \sup_{t\in (0,T)}\iidx{\Og}{\hat{\Fg}^{p_0}}<\infty.\eeq Then together with \mref{interpol} we also have that 
\beqno{interpol1}\itQ{Q}{\hat{\Fg} |v|^2}\le \eg \itQ{Q}{\Fg_0 |Dv|^2}+C(\eg,\bg)\|\Fg\|_{L^1(Q)}\left(\itQ{Q}{ |v|^\bg}\right)^{2/\bg}.\eeq
We can assume that $\hat{\Fg}\ge \Fg$, otherwise the above was proved before. The proof is the same, using the compactness in \refrem{comprem}, under the extra second assumption in \mref{keyFg1}.

\erem

\brem{Fgrem1} Take $\Fg_0=\llg_0>0$. Assume that $v$ satisfies
\beqno{ccFnew}\left|\itQ{Q}{v_t\psi}\right|\le C\itQ{Q}{(|Dv||D\psi|+\Fg|v||\psi|)} \quad \forall\psi\in C^1(Q).\eeq 
With the same conditions on $\Fg$ we also have $$\itQ{Q}{\Fg |v|^2}\le \eg \itQ{Q}{ |Dv|^2}+C(\eg,\bg)\|\Fg\|_{L^1(Q)}\left(\itQ{Q}{ |v|^\bg}\right)^{2/\bg}.$$

\erem

\section{On scalar equations} \label{scalareqn}\eqnoset
We now revisit the regularity theory of scalar equations with integrable coefficients, in the class $\bbM(\Og,T)$. These new improvements  serve well our purposes in the next section.

Note that because $\Og$ is a bounded domain and $T$ is finite, if $\Fg\in\bbM(\Og,T)$, the class of functions defined in \refdef{Ogdef},  then by H\"older's inequality $\Fg$ satisfies the condition  \mref{keyFg} of \reflemm{Interpolation} so that weighted interpolation inequality \mref{interpol} holds. In particular, when $Q=B_R\times(t_0,t_0+R^2)$ we note the following consequence of \mref{interpol}.

\brem{ball} If $Q=B_R\times(t_0,t_0+R^2)$ is a parabolic cylinder ($B_R$ is a ball of radius $R$ in $\Og$) and $\bg=2$. We can make a change of variables $x\to x/R$ and $t\to t/R^2$ to see that \mref{interpol} yields the following inequality
\beqno{interpolloc}\itQ{Q}{\Fg V^2}\le\eg R^{2}\itQ{Q}{\Fg |DV|^2}+C(\eg)R^{-N-2}\|\Fg\|_{L^1(Q)}\itQ{Q}{V^2}.\eeq
\erem

\subsection{Global boundedness and a local estimate}
We consider scalar equation
\beqno{diagSKT} \msysa{v_t=\Div(ADv)+\Div(B_1)+B_2Dv+G_1v+G_2& \mbox{in $Q$,}\\ v=v_0 & \mbox{in $\Og$.}} \eeq
Here $B_1,B_2$ are row vectors in $\RR^N$, $A,G$ are scalar functions.

As usual, we will assume that there is a positive number $\llg_0$ such that
\beqno{ellcond} A\ge \llg_0.\eeq

We also assume that there is a  function $\Fg\in\bbM(\Og,T)$ such that $\Fg\ge\llg_0$ on $Q$ and
\beqno{grate} |A|, |B_1|^2, |B_2|^2, |G_1|, |G_2|\le \Fg.\eeq

By using Steklov average, a weak solution of  \mref{diagSKT} satisfies for all $\eta\in C^1(Q)$
\beqno{weakeqn} \itQ{Q}{u_t\eta+ADvD\eta}=\itQ{Q}{[-B_1D\eta+(B_2Dv+G_1v+G_2)\eta]}.\eeq

We begin with this simple energy lemma.

\blemm{energy} Set $V=|v|^{p}$ with $p\ge1$. For $2p\ge2$, $0<r_1<r_2$ and $0\le t_1<t_2 <T$ let $\fg$ and $\eta(t)$ be cut-off functions respectively in $x,t$ for $B_{r_1}\subset B_{r_2}$ and $(t_1,t_0)\subset (t_2,t_0)$.  Let $Q'=\Og\cap B_{r_1}\times(t_1,t_0)$ and $Q''=\Og\cap B_{r_2}\times(t_2,t_0)$,
 we have
\beqno{keyiter}\barrl{\iidx{B_{r_1}\times\{\tau\}}{V^{2}}+\itQ{Q'}{\Fg |DV|^2}\le C\itQ{Q''}{A |v|^{2p-1}|Dv|\fg|D\fg|}}{2cm} & C\itQ{Q''}{(\Fg V^{2}\fg^2+V^2|D\fg|^2)}+C\itQ{Q''}{|\eta'(t)|V^2}.\earr\eeq
\elemm

\bproof The proof is standard.  We $\eta$ by $|v|^{2p-2}v\fg^2\eta$ in \mref{weakeqn} to get
$$\barrl{\itQ{Q''}{(|v|^{2p}\fg^2\eta)_t}+\itQ{Q''}{A |v|^{2p-2}|Dv|^2\fg^2}\le C\itQ{Q''}{A |v|^{2p-1}|Dv|\fg|D\fg|}}{.5cm} &C\itQ{Q''}{[|B_1||v|^{2p-2}|Dv|\fg^2+|B_1||v|^{2p-1}\fg|D\fg|+|B_2||Dv||v|^{2p-1}\fg^2]}
\\&+\itQ{Q''}{[|G_1||v|^{2p}+|G_2||v|^{2p-1}]\fg^2]}+\itQ{Q''}{|\eta'(t)| |v|^{2p}\fg^2}.\earr$$

Using Young's inequality $$|B_1|v\fg|D\fg|\le |B_1|^2 v^2\fg^2+|D\fg|^2,\; |B_2|v|Dv|\le \eg |B_2| |Dv|^2\fg^2+C(\eg)|B_2|v^2$$ and the growth conditions \mref{grate},  we derive for any $\tau \in (t_2,T)$
\beqno{keyW}\barrl{\iidx{\Og\cap B_{r_1}\times\{\tau\}}{|v|^{2p}}+\itQ{Q}{\Fg |v|^{2p-2}|Dv|^2\fg^2}\le \itQ{Q}{A |v|^{2p}|D\fg|^2}}{2cm} & C\itQ{Q''}{[\Fg |v|^{2p}\fg^2+v^{2p-2}|D\fg|^2]}+C\itQ{Q''}{|\eta'(t)| |v|^{2p}}.\earr\eeq

 We can assume $v\ge1$. The above is \mref{keyiter}. \eproof
 
Applying the usual Moser iteration argument to \mref{keyiter}, we derive
 
\blemm{West} Assume that $\Fg\in\bbM(\Og,T)$ in the growth condition \mref{grate} (and $\Fg\ge\llg_0>0$ on $Q$) and $v$ is a solution of \mref{diagSKT}. Then  there is a constant $C$ depending on $\sup_{(0,T)}\|\Fg\|_{L^{p_0}(\Og)}$ such that \beqno{keyWesta}\sup_{Q}|v|\le C \left(\itQ{Q}{v^2}\right)^\frac12.\eeq

\elemm

\bproof From \reflemm{energy} we prove first that $v$ is globally bounded and so is $A$. Let $r_2>r_1>\mbox{diam}(\Og)$ and $\fg\equiv1$ so that $D\fg=0$. We don't have the first term on the right of \mref{keyiter}.
We choose $\eta$ such that $|\eta'(t)|\le K_0^k$ with $K_0>1$ and $\sum_k \frac{1}{K_0^k}$ converges.
\beqno{eq11}\iidx{\Og\times\{\tau\}}{V^2}+\itQ{Q''}{\Fg |DV|^2}\le C\itQ{Q''}{\Fg V^{2}}+CK_0^k\itQ{Q''}{V^2}.\eeq

\newc{\itQbar}[2]{\dspl{\int\hspace{-2.5mm}\int_{#1}~#2~d\bar{z}}}

Since $v$ satisfies \mref{vccF}, if we multiply the equation of $v$ by $|v|^{2p-2}v\fg$ with $\fg\in C^1(Q)$ then we can see that $V$ also satisfies the equation \mref{vccF} defining $\ccF$ of \reflemm{complem}.

Because $\Fg$ verifies \mref{keyFg}, \reflemm{Interpolation} holds  with $\bg=2$ so that $$\itQ{Q''}{\Fg V^2}\le\eg \itQ{Q''}{\Fg |DV|^2}+C(\eg)\|\Fg\|_{L^1(Q'')}\itQ{Q''}{V^2}.$$ This yields
\beqno{preiter}\iidx{\Og\times\{\tau\}}{V^{2}}+\itQ{Q''}{\Fg|DV|^2}dt\le C(1+K_0^k)\itQ{Q''}{V^{2}}dt.\eeq

Applying the parabolic Sobolev inequality, we get for some $\cg>1$
$$\itQ{Q}{V^{2\cg}}\le C(1+K_0^k)\left(\itQ{Q''}{V^2}\right)^{\cg}.$$

Thus, taking the root, we get $$\left(\itQ{Q}{V^{2\cg}}\right)^\frac{1}{2\cg}\le [C(1+K_0^k)]^\frac{1}{2\cg}\left(\itQ{Q''}{V^2}\right)^\frac12.$$

Using the above inequality with $t_1=t_2-\frac{1}{K_0^k}$ and  the fact that $\sum_k \frac{1}{K_0^k}$ converges, it is standard to apply the Moser iteration technique (let $p$ such that $p=\cg^k$ with $k=0,1,\ldots$ and let $k\to\infty$) to obtain from the above that 
$$\sup_{\Og\times(t_2,T)}|v|\le C^{\frac12\sum_k\cg^{-k}}K_0^{\frac12\sum_k k\cg^{-k}}\left(\itQ{\Og\times(t_1,T)}{|v|^2}\right)^\frac12.$$
The series $\sum_k\cg^{-k},\sum_k k\cg^{-k}$ are convergent so that \mref{keyWesta} holds. \eproof

Similarly, we turn to the local estimate. This type of estimates will be useful for later investigations on the H\"older regularity of weak solutions. We will assume that the function $A$  is bounded. Note that $A$ may depend on $v$ in general and we already showed
that $v$ is bounded globally by the above lemma.

For any $x_0\in\Og$, $R>0$ and $t_0\ge 4R^2$, we define $\Og_R(x_0)=\Og\cap B_R(x_0)$ and $Q_R(x_0)=\Og_R(x_0)\times(t_0-R^2,t_0)$. If $x_0,t_0$ are understood from the context, we simply drop them from the notations.

\blemm{Westloc} Assume that $\Fg\in\bbM(\Og,T)$ in the growth condition \mref{grate} (and $\Fg\ge\llg_0>0$ on $Q$) and that $A$ is bounded. Let $v$ be a solution of \mref{diagSKT} and $B_R$ be a ball in $\RR^N$. Then  there is a constant $C$ depending on $\sup_{(0,T)}\|\Fg\|_{L^{p_0}(\Og)}$ such that \beqno{keyWest}\sup_{Q_R}|v|\le C\left( \frac{1}{R^{N+2}}\itQ{Q_{2R}}{v^2}\right)^\frac12.\eeq

\elemm

In the sequel, we will make use of the local inequality \mref{interpolloc} if we replace $Q$ by $Q_R$ then we need to take a close look at the number $\|\Fg\|_{L^{1}(Q_R)}$. We note that, by H\"older's inequality, $\|\Fg\|_{L^1(\Og_R\times(0,R^2))}\le \sup_t\|\Fg\|_{L^{p_0}(\Og)}R^{2+N(1-1/p_0)}\le \sup_t\|\Fg\|_{L^{p_0}(\Og)}R^{N}$. 

\bproof Let $\fg$ be a cutoff function for $B_R, B_{2R}$. We choose $|D\fg|\le 1/R$ and $\eta$ such that $|\eta'(t)|\le K_0^kR^{-2}$. We have to deal with the first term on the right of \mref{keyiter}. Because $A$ is  bounded, by Young's inequality we have
$$\itQ{Q''}{A |v|^{2p-1}|Dv|\fg|D\fg|}\le \eg\itQ{Q''}{A |v|^{2p-2}|Dv|^2\fg^2}+C(\eg)\itQ{Q''}{|v|^{2p}|D\fg|^2}.$$
Choosing $\eg$ small, the first term on the right hand side can be absorbed into the left.

We now treat the integral of $\Fg V^2\fg^2$. Instead of \reflemm{Interpolation}, \refrem{ball} holds  with $\bg=2$ so that $$\itQ{Q''}{\Fg V^2}\le\eg \itQ{Q''}{\Fg |DV|^2}+C(\eg)\|\Fg\|_{L^1(Q'')}R^{-N-2}\itQ{Q''}{V^2}.$$

Since $\|\Fg\|_{L^1(Q'')}\le \|\Fg\|_{L^{p_0}(\Og)}R^{N}$ (see the note before this proof),  the above inequality, \mref{keyiter} and the choice of $\fg,\eta$ imply a number $C$ depending on $\sup_{(0,T)}\|\Fg\|_{L^{p_0}(\Og)}$ such that  
$$\iidx{\Og_R\times\{\tau\}}{V^{2}}+\itQ{Q''}{\Fg|DV|^2}dt\le C(1+K_0^k)R^{-2}\itQ{Q''}{V^{2}}dt.$$

Making a change of variables $\bar{x}=x/R$, and $\bar{t}=t/R^2$, we see that 
$$R^N\int_{\Og_R\times\{\tau\}}V^{2}d\bar{x}+R^N\itQbar{\Og_R\times(t_2,T)}{\Fg|D_{\bar{x}}V|^2}\le C(1+K_0^k)R^N\itQbar{\Og_{2R}\times(t_1,T)}{V^{2}}.$$
In the new variables $\bar{x}, \bar{t}$, after cancelling $R^N$ in the above, we can assume that $R=1$ in \mref{eq11}. Because $\Fg\ge\llg_0>0$ we get
$$\iidx{\Og_1\times\{\tau\}}{V^2}+\itQ{\Og_1\times(t_2,T)}{|DV|^2}\le C(1+K_0^k)\itQ{\Og_2\times(t_1,T)}{V^2}.$$ 
This is similar to \mref{preiter}. So, we can repeat the iteration argument to obtain
$$\sup_{\Og_1\times(t_2, T)}|v|\le C^{\frac12\sum_k\cg^k}K_0^{\frac12\sum_k k\cg^k}\left(\itQ{\Og_2\times(t_1, T)}{|v|^2}\right)^\frac12.$$
Thus, \mref{keyWest} holds for $R=1$. If we go back to the variables $x,t$ then \mref{keyWest} is also true for any $R>0$.
We finish the proof. \eproof

\subsection{H\"older continuity} 

We will study the H\"older regularity in this subsection. Note that the bounds for the H\"older norm and exponents will depend only on the generic constants and the $L^{p_0}(\Og)$ norms of the parameters in their definition $\bbM(\Og,T)$. This fact will play a crucial role when we estimate the derivatives which appear in cross diffusion systems.

\blemm{blemm} Assume that $\ma$ is bounded and $|\hat{\mb}|^2, |\mb|^2, |\mg|\in \bbM(\Og,T)$. Let $v$ be a weak solution of 
\beqno{eq1}v_t=\Div(\ma Dv)+\hat{\mb}Dv+\Div(\mb)+\mg.\eeq
Then $v$ is H\"older continuous. Its H\"older norm is bounded in terms of the $L^{p_0}(\Og)$ norms of $|\hat{\mb}|^2, |\mb|^2, |\mg|$.
\elemm

The idea based on that of \cite{dleNA}. We present the details and nontrivial modification. Assume first that $\hat{\mb}=\mb=0$. The case $\hat{\mb}\ne0$ and $\mb\ne0$ will be discussed in \refrem{brem} after the proof.

Fixing any $x_0\in\Og$, $t_0>0$ and $4R^2<t_0$, we denote $Q_{iR}=\Og_{iR}\times(t_0-iR^2,t_0)$.

Let $M_i=\sup_{Q_{iR}}v$, $m_i=\inf_{Q_{iR}}v$ and $\og_i=M_i-m_i$. For $\cg>0$ define 
$$ N_1(v)=2(M_4-v)+R^\cg,\quad N_2(v)=2(v-m_4)+R^\cg,$$
$$w_1(v)=\log\left(\frac{\og_4+R^\cg}{2(M_4-v)+R^\cg}\right),\quad w_2(v)=\log\left(\frac{\og_4+R^\cg}{2(v-m_4)+R^\cg}\right).$$

Choosing $\cg$ appropriately, we will prove that either $w_1$ or $w_2$ is bounded from above. This implies a decay estimate for some $\eg\in(0,1)$ and all $R>0$ \beqno{decay}\og_2\le \eg \og_4 + CR^\cg. \eeq It is standard to iterate \mref{decay} to obtain the H\"older continuity of $v$.
Indeed, if either $w_1$ or $w_2$ is bounded from above by $C>0$ in $Q_{2R}$ then this fact implies $$\mbox{either }\og_4+R^\cg\le 2C(\og_4+m_4-v)+2CR^\cg\mbox{ or }\og_4+R^\cg\le 2C(\og_4+v-M_4)+2CR^\cg.$$ 
Taking the supremum (respectively infimum) over $Q_{2R}$ and replacing $m_4$ by $m_2$ (respectively $M_4$ by $M_2$), we obtain
$\og_2\le \eg \og_4+CR^\cg$ for $\eg=\frac{2C-1}{2C}<1$.

Thus, in the sequel, we just need to that either $w_1$ or $w_2$ is bounded from above.

\bproof For any $\eta\in C^1(Q)$ and $\eta\ge0$, observe that $Dw_1=\frac{2Dv}{N_1(v)}$,  $(w_1)_t=\frac{2v_t}{N_1(v)}$
and $Dw_2=-\frac{2Dv}{N_2(v)}$,  $(w_2)_t=-\frac{2v_t}{N_2(v)}$. So, by multiplying the equation of $v$ by $\eta/N_1(v)$ and $-\eta/N_2(v)$ and writing $N_i(v)$ by $N(v)$ we obtain
\beqno{eq2a} \iidx{\Og}{\frac{\partial w}{\partial t}\eta}+ \iidx{\Og}{\myprod{\ma Dw,D\eta}}+\iidx{\Og}{\myprod{\ma Dv,\frac{\eta Dv}{N^2(v)}}}=2\iidx{\Og}{\frac{\pm \mg}{N(v)}\eta}.\eeq

Because $\myprod{\ma Dv, Dv}\ge0$ and $\eta\ge0$,  we can drop the integral of $\myprod{\ma Dv,\eta Dv/N^2(v)}$. We get
\beqno{eq2}\iidx{\Og}{\frac{\partial w}{\partial t}\eta}+ \iidx{\Og}{\myprod{\ma Dw,D\eta}}\le 2\iidx{\Og}{\frac{|\mg|}{N(v)}\eta}.\eeq

Testing \mref{eq2} with $(w^+)^{2p-1}\eta^2$ and using Moser's iteration as in \reflemm{Westloc} with $G_2=\mg/N(v)$. By an appropriate choice of $\cg>0$ (see \mref{mgest} below) we have $\|G_2\|_{L^1(Q_R)}\le C\sup_{(0,T)}\|\mg\|_{L^{p_0}(\Og)}R^N$. So that there is a  constant $C$ depending on $\sup_{(0,T)}\|\mg\|_{L^{p_0}(\Og)}$ such that
\beqno{supw}\sup_{\Og_{2R}\times(t_0-2R^2,t_0)} w^+\le C\left(\frac{1}{R^{N+2}}\itQ{\Og_{4R}\times(t_0-4R^2,t_0)}{w^2}\right)^\frac12.\eeq

If we can show that for any $R>0$ there is a constant $C$ such that
\beqno{wbound} \frac{1}{R^{N+2}}\itQ{\Og_{4R}\times(t_0-4R^2,t_0)}{w^2}\le C\eeq then this implies $w^+$ is bounded. So, the decay estimate \mref{decay} holds. 

Let $\eta$ be a cut-off function for $B_{2R}, B_{4R}$. Replacing $\eta$ in \mref{eq2a} by $\eta^2$, we get
\beqno{k1}\frac{d}{dt}\iidx{\Og}{w\eta^2}+\iidx{\Og}{\ma |Dw|^2\eta^2}\le \iidx{\Og}{\ma|Dw|\eta|D\eta|}+\iidx{\Og}{\frac{|2\mg|}{N(v)}\eta^2}
\eeq

Applying Young's inequality we derive (as $|D\eta|\le C/R$ and $\ma$ is bounded)
\beqno{k2}\barr{lll}\frac{d}{dt}\iidx{\Og}{w\eta^2}+\iidx{\Og}{\ma |Dw|^2\eta^2}&\le& \frac{1}{R^2}\iidx{\Og}{\ma}+\iidx{\Og}{\frac{|2\mg|}{N(v)}\eta^2}\\&\le& CR^{N-2}+\iidx{\Og}{\frac{|2\mg|}{N(v)}\eta^2}.\earr
\eeq

Set $I_*=[t_0-4R^2, t_0-2R^2]$, $Q_*=B_{2R}\times I_*$ and $Q_v=\{(x,t)\in Q_*\,:\,\; w_1\le0\}$. It is easy to see that $w_2\le0$ on $Q_*\setminus Q_v$. Therefore one of $w_1^+,w_2^+$ must vanish on a subset $Q^0$ of $Q_*$ with $|Q^0|\ge \frac12|Q_*|$. We denote by $w$ such function. Let $Q_t^0$ be the slice $Q^0\cap (B_{2R}\times\{t\})$ then $Q^0=\cup_{t\in I_*} Q_t^0$.  For $t\in I_*$ let
$$\Og^0_t=\{x\,:\, w^+(x,t)=0\}, \quad m(t)=|Q^0_t|.$$
The fact that $|Q^0|\ge \frac12|Q_*|$ implies $\int_{I_*}m(t)dt\ge \frac12 R^{N+2}$.

We now set $$V(t)=\frac{\iidx{\Og}{w\eta^2}}{\iidx{\Og}{\eta^2}}.$$

By the weighted Poincar\'e' inequality (\cite[Lemma 3]{Moser}) 
$$\iidx{\Og}{(w-V)^2\eta^2}\le CR^2\iidx{\Og}{|Dw|^2\eta^2}.$$
Reducing the integral on the left to the set $\Og^0_t$ where $w\le0$ (so that $V^2\le (w-V)^2$), we have $$V^2(t)m(t)\le CR^2\iidx{\Og}{|Dw|^2\eta^2}.$$
Since $N(v)\ge R^\cg$ on $Q_{4R}$, the above estimate and \mref{k2} implies that ($V'$ denotes the $t$ derivative)
\beqno{k3}R^N V'(t)+\frac{1}{R^2}V^2(t)m(t)\le CR^{N-2}+\frac{2}{R^\cg}\iidx{\Og}{|\mg|\eta^2}.\eeq

Because of $\mg\in\bbM(\Og,T)$ for some $q>N/2$, by an appropriate choice of $\cg>0$, we also get
\beqno{mgest}\iidx{\Og}{|\mg|\eta^2}\le R^{N(1-\frac{1}{q})}\left(\iidx{\Og}{|\mg|^q}\right)^\frac{1}{q}\le CR^{N-2+\cg}.\eeq 
This also shows that $\|\mg/N(v)\|_{L^1(Q_R)}$ is bounded by $CR^N$.

Thus, from \mref{k3}
\beqno{k4}R^N V'(t)+\frac{1}{R^2}V^2(t)m(t)\le CR^{N-2}.\eeq

We show that $V(t_1)$ is bounded on for some $t_1\in I_*$. Indeed, supppose $V(t)\ge A>0$ in $I_*$. We have from \mref{k3} and \mref{mgest}
$$ R^{N+2}\frac{V'(t)}{V^2(t)}+m(t)\le \frac{CR^{N}}{A^2}.$$

Because $\int_{I_*}m(t)dt\ge \frac12 R^{N+2}$ and $|I_*|\sim R^2$, we integrate this over $I_*$ to see that
$$R^{N+2}\le \int_{I_*}m(t)dt\le R^{N+2}\left(\frac{2}{A}+\frac{C}{A^2}\right).$$
By choosing $A$ large we get a contradiction. So, we must have $V(t_1)\le A$ for some $t_1\in I_*$.

Integrating \mref{k2} over $[t_1,t_2]$ for any $t_2\in I_0=[t_0-2R^2, t_0]$, we have
$$V(t_2)\iidx{\Og\times\{t_2\}}{\eta^2}+\itQ{\Og\times I_0}{\ma |Dw|^2\eta^2}\le CR^N +V(t_1)\iidx{\Og\times\{t_1\}}{\eta^2}.$$
This implies that $V(t_2)\le C$ for all $t_2\in I_0$ and $\itQ{\Og\times I_0}{\ma |Dw^+|^2\eta^2}\le CR^N$. This implies that $V(t_2)\le C$ for all $t_2\in I_0$ and $\itQ{\Og\times I_0}{\ma |Dw|^2\eta^2}\le CR^N$. Since we can always assume that $\og_4\ge R^\cg$ (otherwise there is nothing to prove) to have that $w\ge\log2$. Thus, $V(t)$ is bounded from below so that $|V(t)|$ is bounded.

By Poincar\'e's inequality again we have
$$\itQ{\Og\times I_0}{(w-V)^2\eta^2}\le CR^2\itQ{\Og\times I_0}{|Dw|^2\eta^2}\le CR^{N+2}.$$

Since $|V(t)|$ is bounded on $I_0$, replacing $R$ by $2R$, the above implies the desired \mref{wbound}. \eproof

\brem{brem} Assume that $\mb\ne0$. Replacing $\eta$ by $\eta/N(v)$ we have the extra terms in \mref{eq2a}
$$\iidx{\Og}{\myprod{\mb, \frac{D\eta}{N(v)}+\frac{\eta Dv}{N^2(v)}}}=\iidx{\Og}{\myprod{\frac{\mb}{N(v)}, D\eta+\frac12 Dw \eta}}.$$

We take $\eta$ to be $(w^+)^{2p-2}w^+\eta^2$ and think of $B_1:= \mb/N(v), B_2:=\frac12\mb/N(v)$ (see \mref{mbest} below). Again noting that $Dw=2Dv/N(v)$. Because $|\mb|^2 \in \bbM(\Og,T)$, the Moser's technique is applied as before. We  get the local estimate \mref{supw}.

When we estimate $R^{-2N-2}\itQ{Q_{2R}}{w^2}$ to show that $w$ is bounded, for $\eta$ is a cut-off function with $|D\eta|\le 1/R$ we need the following estimate for the extra terms in \mref{k1} \beqno{best}\iidx{\Og}{\frac{\mb D\eta}{N(v)}},\; \iidx{\Og}{\frac{\mb^2 \eta}{N^2(v)}}\le R^{N-2}.\eeq

The second term appears due to an use of Young's inequality.

Because $N(v)\ge R^\cg$, by H\"older's inequality we have $$\iidx{\Og}{\frac{\mb D\eta}{N(v)}}\le \frac{1}{R^{1+\cg}}\iidx{\Og_{2R}}{\mb}\le CR^{-1-\cg}(R^N)^\frac12\left(\iidx{\Og_{2R}}{\mb^2}\right)^\frac12$$
Another use of H\"older's inequality again and the assumption that $|\mb|^2\in \bbM(\Og,T)$ for some $p_0>N/2$ the last integral is bounded by $\||\mb|^2\|_{L^{p_0}(\Og)}R^{N-2+\cg_0}$ for some $\cg_0>0$. So that
$$\iidx{\Og}{\frac{\mb D\eta}{N(v)}}\le CR^{-1-\cg+\frac{N}{2}+\frac{N-2+\cg_0}{2}}=CR^{N-2+\frac{\cg_0}{2}-\cg}.$$
By the same reason 
\beqno{mbest}\iidx{\Og}{\frac{\mb^2 }{N^2(v)}\eta}\le CR^{-2\cg+N-2+\cg_0}.\eeq

An appropriate choice of choice of $\cg>0$ we obtain \mref{best}. This also shows that $\|\mb^2/N^2(v)\|_{L^1(Q_R)}\le CR^N$ for a constant depending on $\sup_{(0,T)}\|\mb^2\|_{L^{p_0}(\Og)}$ (compare with the proof for $\mg$ in \mref{mgest}) and again confirms that the constant $C$ in \mref{supw} is depending on $\sup_{(0,T)}\|\mb^2\|_{L^{p_0}(\Og)}$.

If $\hat{\mb}\ne0$ then we also have the extra term 
$$\iidx{\Og}{\myprod{\hat{\mb}Dv, \frac{\eta}{N(v)}}}=\iidx{\Og}{\myprod{\frac12\hat{\mb}Dw,  \eta}}$$ in \mref{eq2a}. The argument is similar and much simpler.

\erem

\subsection{Estimates for derivatives} 

\blemm{Dup} If $v$ is H\"older continuous then for some $R_0$ small and any $p$ such that $2p>1$ and $0<t_0<t<T$ 
\beqno{Dupest1}\sup_{(t_0,t)}\iidx{\Og}{|Dv|^{2p}}\le  C(R_0^{-2},t_0^{-1})\itQ{\Og\times(0,t)}{|Dv|^{2p}}.\eeq

Moreover, \beqno{Dupest2}\sup_{(t_0,t)}\iidx{\Og}{|Dv|^{2p}}\le  C\left(R_0^{-2},t_0^{-1},\itQ{\Og\times(0,t)}{|Dv|^{2}}\right).\eeq

\elemm

\bproof Differentiate the equation of $v$ to see that \beqno{Dueq1} (Dv)_t=\Div(\ma D^2v+\ma_v\myprod{Dv,Dv}+\hat{\mb}_v\myprod{Dv,Dv}+\hat{\mb}D^2v+\Div(\mb_vDv)+\mg_v  Dv.\eeq

Since the parameters $\ma,\hat{\mb}, \mb$ and $\mg$  of the equation are bounded, we test the equation with $|Dv|^{2p-2}Dv\fg^2$ and use Young's inequality to get
\beqno{Dueq2}\barrl{\iidx{\Og}{|Dv|^{2p}\fg^2}+\itQ{Q}{|Dv|^{2p-2}|D^2v|^2\fg^2}\le}{3cm}& C\itQ{Q}{|Dv|^{2p+2}\fg^2}+C\itQ{Q}{|Dv|^{2p}(\fg^2+|D\fg|^2)}.\earr\eeq
The key issue here is that we will have to handle  the integral of $|Dv|^{2p+2}$.

For any $t>0$ and $0<t_0<T$ we write $\Og_t=\Og\cap B_t$, $Q_t=\Og_t\times(t_0,T)$ $$\ccA_p(t)=\sup_{(t_0,T)}\iidx{\Og_t}{|Dv|^{2p}},\;\ccB_p(t)=\itQ{Q_t}{|Dv|^{2p+2}},$$
$$\ccH_p(t)=\itQ{Q_t}{|Dv|^{2p-2}|D^2v|^2},\; \ccG_p(t)=\itQ{Q_t}{|Dv|^{2p}}.$$

By the local Gagliardo-Nirenberg inequality (\cite[Lemma 2.4]{dleANS} with $u=U$) we showed that $$\iidx{\Og_s}{|Dv|^{2p+2}}\le C\|v\|_{BMO(\Og_t)}^2\iidx{\Og_t}{|Dv|^{2p-2}|D^2v|^2}+C\frac{\|v\|_{BMO(\Og_t)}}{(t-s)^2}\iidx{\Og_t}{|Dv|^{2p}}.$$

Using this inequality in \mref{Dueq2}, where $\fg$ is a cut-off function for $B_s,B_t$ and $(t_0,T)$, by H\"older continuity $\|v\|_{BMO(\Og_t)}$ can be very small. From \mref{Dueq2}, for $0<s<t<R_0$ with sufficiently small $R_0$ there is $C=C(t_0^{-1})$ the following recursive system of inequalities hold
$$\ccA_p(s)+\ccH_p(s)\le C\ccB_p(t)+ \frac{C}{(t-s)^2}\ccG_p(t),$$
$$\ccB_p(s)\le \eg \ccH_p(t)+\frac{C}{(t-s)^2}\ccG_p(t).$$
If  $\eg$ is small (that is $R_0$ is small because $v$ H\"older continuous) then we can iterate this to absorb the terms $\ccH_p,\ccB_p$ on the right hand side to the left hand side to get (see \cite[inequality (3.27), proof of Proposition 3.1]{dleANS})
\beqno{Duiter}\ccA_p(R_0)+\ccH_p(R_0)\le \frac{C}{R_0^2}\ccG_p(2R_0).\eeq

This gives the estimate of $\ccA_p$. That is
$$\sup_{(t_0,T)}\iidx{\Og_{R_0}}{|Dv|^{2p}}\le CR_0^{-2}\itQ{Q_{2R_0}}{|Dv|^{2p}}.$$

Fixing such $R_0$ and partitioning $\Og$ into balls of radius $R_0$ and summing the above inequality over this partition, we derive a global estimate
\beqno{DuLp}\sup_{(t_0,T)}\iidx{\Og}{|Dv|^{2p}}\le C(R_0^{-2},t_0^{-1})\itQ{\Og\times(t_0,T)}{|Dv|^{2p}}.\eeq
This gives the first assertion \mref{Dupest1} of the lemma.

By the parabolic Sobolev inequality, \mref{Duiter} also shows that its right hand side is self-improved. We recall the parabolic Sobolev inequality
$$\itQ{\Og\times(0,T)}{|g|^\frac{2}{N}|G|^2} \le C\sup_{(0,T)}\left(\iidx{\Og}{|g|}\right)^{\frac{2}{N}}\left(\itQ{\Og\times(0,T)}{|DG|^2+|G|^2}\right).$$
which can be proved easily by H\"older's and Sobolev's inequalities. Let $g=G=|Dv|^p$. We see that
 if $\ccG_p(2R_0)$ is finite for some $p\ge1$ then so are $\ccA_p(R_0), \ccH_p(R_0)$ so that $\ccG_{\cg p}(R_0)$ is also finite for  some $\cg=1+\frac{2}{N}>1$. We can repeat this argument (finite times) to see that $\ccG_p(R_0)$ is finite for any $p>1$.

Testing the equation of $v$ by $v$, we see easily that the right hand side of \mref{DuLp} is finite for $p=1$. The second assertion \mref{Dupest2} then follows. \eproof

\brem{Duvrem} Because $v$ is bounded, the factors in \mref{Dueq1} are all bounded. We can allow to include  terms $\mF$, $\hat{\mF}Dv$, with integrable $\mF,\hat{\mF}$, in \mref{Dueq1}. In this case, applying  Young's inequality to the term $\mF|Dv|^{2p-1}$, with exponents $\frac{2p+2}{2p-1}$ and $\frac{2p+2}{3}$ ($\frac{2p+2}{2p}$, $\frac{2p+2}{2}$ for $\hat{\mF}|Dv|^{2p}$), \mref{Dueq2} becomes $$\barrl{\iidx{\Og}{|Dv|^{2p}\fg^2}+\itQ{Q}{|Dv|^{2p-2}|D^2v|^2\fg^2}\le C\itQ{Q}{|Dv|^{2p+2}\fg^2}+}{2cm}& C\itQ{Q}{|f_v||Dv|^{2p}\fg^2} + \itQ{Q}{(|\mF|^\frac{2p+2}{3}+|\hat{\mF}|^\frac{2p+2}{2})\fg^2}.\earr$$

By the same argument, we obtain another version of \mref{Duiter} (see \cite[Lemma B.1.7]{dlebook})
$$\ccA_p(R_0)+\ccH_p(R_0)\le \frac{C}{R_0^2}\ccG_p(2R_0)+\itQ{\Og_{2R_0}\times(t_0,T)}{(|\mF|^\frac{2p+2}{3}+|\hat{\mF|}|^\frac{2p+2}{2})},$$ with $\itQ{Q_t}{(|\mF|^\frac{2p+2}{3}+|\hat{\mF}|^\frac{2p+2}{2})}$ playing the role of $h$ in the equations (B.1.11) and (B.1.12) of \cite{dlebook}) 
This shows that  
$$\barr{lll}\sup_{(t_0,T)}\iidx{\Og}{|Dv|^{2p}}&\le& C(R_0^{-2},t_0^{-1})\itQ{\Og\times(t_0,T)}{|Dv|^{2p}}+\\&&\itQ{\Og\times(t_0,T)}{(|\mF|^\frac{2p+2}{3}+|\hat{\mF}|^\frac{2p+2}{2})}.\earr$$

As before, we note that the quantity $$\itQ{\Og\times(t_0,T)}{|Dv|^{2p}}$$ is still self-improved so that the estimate for $\sup_t\ccA_{2p}(R_0)$ holds as long as  $$\itQ{\Og\times(t_0,T)}{(|\mF|^\frac{2p+2}{3}+|\hat{\mF}|^\frac{2p+2}{2})}<\infty.$$

\erem

\section{Applications to systems}\label{sys}{\eqnoset}

In this section, we apply the theory of scalar equations to systems of the form
\beqno{fullsys0}\msysa{W_t=\Div(\mA DW)+\mB DW+\mG W+F&\mbox{in $Q$,}\\ \mbox{Homogeneous Dirichlet or Neumann boundary conditions}&\mbox{on $\partial\Og\times(0,T)$},
\\
W=W_0&\mbox{on $\Og$.}}\eeq
where $W=[u,v]^T$, $\mA=\mat{\ag&\bg\\\cg&\dg}$, $\mB=\mat{a&b\\c&d}$, $\mG=\mat{g_1&0\\0&g_2}$ and $F=\mat{f_1\\f_2}$. The entries of these matrices are functions in $u,v$.

As usual, we always assume the ellipticity condition that there are some function $\llg$ and a positive constant $\llg_0$ such that $\llg\ge\llg_0$ and
\beqno{ellcond}\myprod{\mA DW,DW}\ge \llg(W)|DW|^2 \quad \forall W\in C^1(\Og,\RR^2).\eeq

The global existence problem of the system \mref{fullsys0} can be established by Amann's theory in \cite{Am2}. We just need to show that \beqno{est1} \limsup_{t\to T}\|W\|_{W^{1,2p}(\Og)} <\infty \mbox{ for any $T$ and some $p>N/2$.}\eeq

Let $\llg\le \LLg$ be eigenvalues of $\mA$ and consider  the ratio $\nu=\frac{\llg}{\LLg}$. We proved in \cite{dletrans, dleANS} that if $s>-1$ and $\frac{s}{s+2}<\nu$ then \mref{ellcond} implies a $c>0$ such that
\beqno{ellcond1}  \myprod{\mA DX,D(|X|^{s}X)}\ge c\llg(X)|X|^{s}|DX|^2, \quad \forall X\in C^1(\Og,\RR^2).\eeq

{\bf Example: } For simplicity, assume that the entries of $\mA,\mB, \mG$ are linear functions of $u,v$, $F\equiv0$ and $\llg(W)\sim \llg_0+|W|$. Assume also $\frac{s}{s+2}<\nu$ so that \mref{ellcond1}  is available.

Using the above assumption, if $1-\frac{2}{2p-1}<\nu$ we test \mref{fullsys0} with $|W|^{2p-3}W$ ($p>1$) and use \mref{ellcond1} with $X=W$ and $s=2p-3$ and the fact that $\llg(W)\sim |W|$ to get the following energy estimate
$$\iidx{\Og}{|W|^{2p}}+\itQ{Q}{\llg(W)|W|^{2p-2}|DW|^2}\le \iidx{\Og}{|W_0|^{2p}} +C_0\itQ{Q}{|W|^{2p}}+C_1.$$
Here, the constant $C_0$ depends on the coefficients of $u,v$ of the entries of $\mB,\mG$.

From the Sobolev inequality for $|W|^{p}$ we have a universal constant $c_0$ such that
$$c_0\iidx{\Og}{|W|^{2p}}\le \iidx{\Og}{|W|^{2p-2}|DW|^2}.$$

If $C_0<c_0$ (this is true if the coefficients of the entries of $\mB,\mG$ are small in comparing with $\llg_0$) then 
$$\iidx{\Og}{|W|^{2p}}\le \iidx{\Og}{|W_0|^{2p}}+C_3.$$

Together with the energy estimate this implies a Gr\"onwall inequality for $\|W\|_{L^{2p}(\Og)}$ and provides a uniform estimate for it if $W_0\in L^{2p}(\Og)$. Thus, if we assume the condition $1-\frac{2}{N-1}<\nu$ (of course this is void if $N\le3$ because $\nu\in (0,1)$) then we can take $p>N/2$ and the conditions on the membership in $\bbM(\Og,T)$ of the parameters of the systems are verified so that the theory for scalar equations in the previous section can be applied to each equation. 

The above argument can be extended to the case when the matrices of \mref{fullsys0} are of any sizes. In this case, we will denote by $\llg, \LLg$ to be the smallest and largest eigenvalues of $\mA$.

{\bf On the triangular systems:} 

We consider the system \beqno{apptri}W_t=\Div(\mA DW)+\mB DW+\mG W+F ,\eeq
where $\mA, \mB$ are tringular matrices. That is, $\mA=\mat{\ag&\bg\\0&\dg}$, $\mB=\mat{a&b\\0&d}$.

We will write the equations of $u,v$ as scalar ones of the forms \mref{diagSKT} or \mref{eq1} discussed in the previous section
$$v_t=\Div(ADv)+\Div(B_1)+B_2Dv+G_1v +G_2,$$
$$v_t=\Div(\ma Dv) +\hat{\mb}Dv+\Div(\mb)+\mg.$$

Let us summarize our findings there for the readers convenience in the proofs below. \reflemm{West} shows that if the data $A, B_i^2, G_i$ of \mref{diagSKT} are in the class $\bbM(\Og,T)$ then its solutions are bounded globally on $(0,T)$. The H\"older continuity of the solutions of \mref{eq1}  is guaranteed by \reflemm{blemm} if $\hat{\mb}^2, \mb^2, \mg\in\bbM(\Og,T)$. The integrability their derivatives over $\Og$ is investigated in \reflemm{Dup} and its \refrem{Duvrem}.

In this section, we will always assume that {\em $\dg$ depends only on $v$} and
\beqno{tri1} \dg, d^2, g_2, f_2\in \bbM(\Og,T),\eeq
\beqno{tri2} \ag, a^2, g_1, f_1\in \bbM(\Og,T).\eeq
Also, their $L^{p_0}$ norms in the class $\bbM(\Og)$ are uniformly bounded.

We will begin with the simplest case when the parameters $d,g_2, f_2$ of the equation of $v$ depend only on $v$.

\btheo{trithm1} Assume \mref{tri1}, \mref{tri2} and that $d,g_2, f_2$ depend only on $v$. If $\bg^2, b\in\bbM(\Og,T)$ then the problem \mref{apptri} has a unique strong solution on $(0,T)$ for any $N$. \etheo

\bproof
Since the parameters of the $v$ equation depends only on $v$ and \mref{tri2}, from \reflemm{West} and \reflemm{blemm} we see that $v$ is bounded and H\"older continuous. Also, by \reflemm{Dup} $Dv\in L^p(\Og)$ for all $p>1$. From the equation of $u$, let $B_1:=\bg Dv$, $B_2=a$, $G_1=g_1$ and $G_2=bDv+f_1$. The assumption that $\bg^2,b\in \bbM(\Og,T)$ then implies $B_1^2, B_2^2, G_1,G_2^2\in\bbM(\Og,T)$ for another $p_0>N/2$. So that $u$ is bounded by \reflemm{West} is H\"older continuous by \reflemm{blemm}. We conclude that $u,v$ exist on $(0,T)$. \eproof

On the other hand, if we allow $d, g_2,f_2$ to depend on $u$ and somehow can establish suitable integrabilities of $|d_u||Du|,|(g_2)_u||Du|, |(f_2)_u||Du|$. Then there will be  extra terms $\hat{\mF}Dv, \mF$ on the right hand side of \mref{Dueq1} (when we differentiate the scalar equation in $v$). Under some suitable integrabilty assumptions, this implies that $|Dv|\in L^{2p}(\Og)$ by \refrem{Duvrem}.

We have the following general result when $d, g_2, f_2$ depend on $u,v$. 
\btheo{genN} Suppose  the spectral gap condition such that \mref{ellcond1} holds. Assume also that for some $p>N/2$ \beqno{triFa}|d_u Du|\in  L^\frac{2p+2}{2}(Q)\mbox{ and }|(g_2)_uDu|,\;|(f_2)_uDu|\in L^\frac{2p+2}{3}(Q).\eeq

Assume also that there is $q$ such that $p>q>N/2$ and
\beqno{tribg} \iidx{\Og}{|\bg|^\frac{2qp}{p-q}},\;\iidx{\Og}{|b|^\frac{2qp}{2p-q}}< \infty.\eeq

Then the problem \mref{apptri} has a unique strong solution  on $(0,T)$ for any $N$. \etheo

\bproof The assumption \mref{tri1} implies that $v$ is bounded and H\"older continuous. In applying \reflemm{Dup}, the dependence of $d,g_2,f_2$ gives rise to some extra terms in \mref{Dueq1}. We need to refer to \refrem{Duvrem} in what follows.
 \refrem{Duvrem} and \mref{triFa}, with $\hat{\mF}=d_uDu$ and $\mF= (g_2)_uDuv+(f_2)_uDu$, yield an estimate for $\|Dv\|_{L^{2p}(\Og)}$ for some $p>N/2$.

The equation of $u$ can be written as a scalar equation of the form \mref{diagSKT} with $B_1:=\bg Dv$, $B_2=a$, $G_1=g_1$ and $G_2=bDv+f_1$.  If $p>q$ then H\"older's inequality gives
$$\iidx{\Og}{|B_1|^{2q}}\le \left(\iidx{\Og}{|\bg|^\frac{2qp}{p-q}}\right)^{1-\frac{q}{p}}\left(\iidx{\Og}{|Dv|^{2p}}\right)^\frac{q}{p},$$
$$\iidx{\Og}{|bDv|^{q}}\le \left(\iidx{\Og}{|b|^\frac{2qp}{2p-q}}\right)^{1-\frac{q}{2p}}\left(\iidx{\Og}{|Dv|^{2p}}\right)^\frac{q}{2p}.$$

By \mref{tribg} $B_1^2, G_2\in\bbM(\Og,T)$, as $q>N/2$. Together with \mref{tri2},  this imply that $u$ is bounded.

The equation of $u$ can also be written as a scalar equation of the form \mref{eq1} with $\mb:=\bg Dv$. From \mref{tribg} we also have $\mb^2\in\bbM(\Og,T)$ ($u,v$ are now bounded) \reflemm{blemm} gives that $u$ is H\"older continuous.

Now, as $u,v$ are bounded and H\"older continuous, if the spectral gap condition such that \mref{ellcond1} holds then the argument in \reflemm{Dup}  also gives an estimate for $\|Du\|_{L^{2p}(\Og)}$. The estimate for $\|DW\|_{L^{2p}(\Og)}$ with $W=[u,v]^T$ and $2p>N$ then completes the proof. \eproof

In particular, if $N=2,3$ then the condition \mref{triFa} is  guaranteed. We easily have

\btheo{23thm} Assume that  there are $q_1, q_2$ such that $2\ge q_1>\frac{N+2}{2}$ and $2\ge q_2>\frac{N+2}{3}$ such that ($\llg$ is the ellipticity function)
\beqno{dgfcond}|d_u|^{q_1},\;|(g_2)_u|^{q_2},\; |(f_2)_u|^{q_2}\le \llg.\eeq
Here, $\llg$ is the ellipticity function in \mref{ellcond}. Then \mref{triFa} holds for some $p>N/2$.

Suppose also that \mref{tribg} holds for some $p>q>N/2$. 
The problem \mref{apptri} then has a  strong solution on $(0,T)$. \etheo

\bproof  Testing the equations of $u,v$ by $u$ and $v$ respectively and using \mref{ellcond}, we easily see that $ \llg|Du|^2\in L^1(Q)$. 

We define $\hat{\mF}=|d_u||Du|$ and $\mF= |(g_2)_u||Du|+|(f_2)_u||Du|$. It is easy to see that  then we can find $p>N/2$ such that  $q_1=\frac{2p+2}{2}$ and $q_2=\frac{2p+2}{3}$. Clearly, the condition \mref{dgfcond} implies that $\hat{\mF}^{q_1}, \mF^{q_2}\le \llg|Du|^2$ if $q_1,q_2\le  2$. \reflemm{Dup} and its \refrem{Duvrem} then yield an estimate for $\|Dv\|_{L^{2p}(\Og)}$.

We see that one can find $q_1,q_2\le 2$ if $N=2$. If $d$ is independent of $v$ then we just need $q_2\le 2$. This is the case if $N=3$.

We see that our theorem follows from \reftheo{genN}. \eproof

\brem{drem} One of the main reason we have to assume that $\mB$ is upper triangular is that in order to use \refrem{Duvrem} to obtain an estimate of $Dv$ we have to differentiate the equation of $v$. In doing so, if $\mB$ is upper triangular then we do not have a term involving $D^2u$, whose integrability is not known, so that the argument can go on under the assumption \mref{triFa}. Similarly, we have to assume that $\mA$ is upper triangular and $\dg$ must be independent of $u$ as we don't know yet $Du\in L^{2p+2}(\Og)$, a fact we would need in \refrem{Duvrem} if $\dg$ had depended on $u$ (see \refrem{deltauv}).

\erem

\brem{deltauv} We can allow $\dg$ to depend on both $u,v$ as long as it is bounded in terms of $u\in\RR$. As before, if $\dg, f_2\in \bbM(\Og,T)$ then $v$ is globally bounded by \reflemm{West} and then H\"older continuous by \reflemm{blemm}. The key point is that we need an estimate for $\|Dv\|_{L^{2p}(\Og)}$.  From the proof of \reflemm{Dup} we have an extra term $\ma_uDu Dv$ (which is now $\dg_u DuDv$) on the right of \mref{Dueq2}. As $\ma_u$ is now bounded, after testing by $|Dv|^{2p-2}Dv\fg^2$, by Young's inequality we have the term $|Dv|^{2p}|Du|^2$. Again, we have  $ |Dv|^{2p}|Du|^2\le C(|Dv|^{2p+2}+|Du|^{2p+2})$. The integral of $|Dv|^{2p+2}$ will be treated as before. We will make an assumption $$\itQ{Q}{|Du|^{2p+2}}<\infty \mbox{ for $p>N/2$.}$$
Then we still have the needed estimate for $\sup_{(0,T)}\|Dv\|_{L^{2p}(\Og)}$ as before to use in the equation for $u$.

\erem

\reftheo{genN} is easily extended to triangular systems of more than two equations. We consider the system
\beqno{trisys} u_t=\Div(\mA Du)+\mB Du +\mG u.\eeq
Here, $\mA=[\ag_{ij}]$, $\mB=[\bg_{ij}]$ are upper triangular matrices ($\ag_{ij}=\bg_{ij}=0$ if $i<j$) and $\mG=[g_{ij}]$ is a full matrix and $F=f_i$. 
In particular, we will assume that $\ag_{ij}$ depends only on $u_j$ for $j\ge i$ (see \refrem{drem}).

We restate \reftheo{trisysthma} in the Introduction here and present its proof.
\btheo{trisysthm} Assume the following integrability conditions
\beqno{trisys1} \ag_{i1}, \bg_{i,1}^2, f_{i}\in \bbM(\Og,T).\eeq

Suppose that for some $p>N/2$ and any $i<j$ \beqno{triFsys}|(\bg_{ij})_{u_i}Du_i|\in L^\frac{2p+2}{2}(Q) \mbox{ and }|(g_{ij})_{u_i}Du_i|,\;|(f_j)_{u_i}Du_i|\in L^\frac{2p+2}{3}(Q).\eeq

And for some $q$ such that $p>q>N/2$ and any $i<j$
\beqno{trisysbg} \iidx{\Og}{\ag_{i,j}^\frac{2qp}{p-q}},\; \iidx{\Og}{\bg_{i,j}^\frac{2qp}{2p-q}}< \infty, \eeq

Then \mref{trisys} has a unique strong solution on $(0,T)$ for any $N$.
\etheo

\bproof By induction (backward), suppose that we have shown that $u_i$ is bounded and H\"older continuous for $i\ge k$. This is true for $i=k=m$ by \reftheo{genN}.
By \mref{triFsys} and \refrem{Duvrem}, with $\hat{\mF}=\sum_{i<j}|(\bg_{ij})_{u_i}Du_i|$ and $\mF:=\sum_{i<j}|(g_{ij})_{u_i}Du_i|+|(f_{j})_{u_i}Du_i|$, yields that $|Du_j|\in L^{2p}(\Og)$ for some $p>N/2$ and any $j\ge k$.

The equation of $u_{k-1}$ can be written as a scalar equation of the form \mref{diagSKT} with  $A:=\ag_{k-1,1}$, $B_1:=\sum_{i\ge k}\ag_{k-1,i} Du_i$, $B_2:=\bg_{k-1,k-1}$ and $G_2:=\sum_{i\ge k}\bg_{k-1,i} Du_i$, where $i\ge k$.  If $p>q$ then H\"older's inequality gives
$$\iidx{\Og}{B_1^{2q}}\le \sum_{i\ge k}\left(\iidx{\Og}{\ag_{k-1,i}}^\frac{2qp}{p-q}\right)^{1-\frac{q}{p}}\left(\iidx{\Og}{|Du_i|^{2p}}\right)^\frac{q}{p},$$
$$\iidx{\Og}{|\bg_{k-1,i}Du_i|^{q}}\le \left(\iidx{\Og}{\bg_{k-1,i}^\frac{2qp}{2p-q}}\right)^{1-\frac{q}{2p}}\left(\iidx{\Og}{|Du_i|^{2p}}\right)^\frac{q}{2p} \quad \forall i\ge k.$$

By  \mref{trisysbg}  $B_1^2, B_2^2, G_2\in\bbM(\Og,T)$, as $q>N/2$. Together with \mref{trisys1}, $A\in \bbM(\Og,T)$,  we can conclude that $u_{k-1}$ is bounded.

The equation of $u_{k-1}$ can also be written as a scalar equation of the form \mref{eq1} with $\hat{\mb}=\bg_{k-1,1}$, $\mb:=\sum_{i\ge k}\ag_{k-1,i} Du_i$. The above estimate of $B_1$ also yields $\mb^2\in\bbM(\Og,T)$. \reflemm{blemm} then gives that $u_{k-1}$ is H\"older continuous. Thus, $u_i$ is bounded and H\"older continuous for all $i$. We also have $\|Du_i\|_{L^{2p}(\Og)}$ is bounded for some $p>N/2$. 
The proof is complete. \eproof

The similar version of \reftheo{23thm} for large systems is

\btheo{23trithm} Assume that  $|(f_i)_{u_j}|^2\le \ag_{i-1,1}$ for $1\le j< i$. Suppose that $N=2,3$ and $\bg\in L^{q_0}(\Og)$ for all $t\in (0,T)$  and for some $q_0>\frac{4N}{4-N}$. Then the problem \mref{apptri} has a unique strong solution on $(0,T)$ for any $N$. \etheo

{\bf On a special full system:} Let $\mA$ be a full matrix
$$\mA=\mat{a&b\\c &d}.$$
We want to find the matrices $\ma=\mat{\ag&\bg\\0&\dg}$ and $\ccJ=\mat{s_1&s_2\\s_3&s4}$ such that $\mA=\ccJ^{-1}\ma\ccJ$. In this case, $\ccJ\mA\ccJ^{-1}=\ma$, a triangular matrix.

Let $D=\mbox{det}\ccJ$. Straightforward calculations show that $\ccJ\mA\ccJ^{-1}=\ma$ is equivalent to
\beqno{abcd}\left\{\barr{c}(s_1a+s_2c)s_4-(s_1b+s_2d)s_3=\ag D,\\-(s_1a+s_2c)s_2+(s_1b+s_2d)s_1=\bg D,\\(s_3a+s_4c)s_4-(s_3b+s_4d)s_3=0,\\
-(s_3a+s_4c)s_2+(s_3b+s_4d)s_1=\dg D.
\earr\right.\eeq

The third and fourth equations  are equivalent to $$\left\{\barr{c}s_3a+s_4c=s_3\dg,\\s_3b+s_4d=s_4\dg.\earr\right.$$ 

Suppose that there is a number $k\ne0$ such that \beqno{s3s4}\dg:=a+kc=\frac{1}{k}b+d\eeq
is a positive function in $v$. We then choose $s_3,s_4$ such that $s_4/s_3=k$. The above system is verified.

For such $s_3,s_4$ we can choose $s_1, s_2$ and define $\ag,\bg$ from the first two equations of \mref{abcd} such that $\ag>0$. We then see that $V=\ccJ W$ satisfies a system of the
form $V_t=\Div(\ma DV)+\hat{\mg} V$ with $\ma$ is triangular. The previous argument shows that $V$ H\"older continuous (so that $v$ exists globally). This applies also to $W$.

{\bf Proof of \reftheo{fullthma}:}  We assume only that $\sup_{(0,T)}\|Du\|_{L^{2p}(\Og)}<\infty$ for some $p>N/2$ (so that $u$ is bounded). Since we also suppose that for some $q$ such that $p>q>N/2$ and
$$ \iidx{\Og}{\cg^\frac{2qp}{p-q}}  <\infty,$$
we can write the equation of $v$ in the form \mref{diagSKT} or \mref{eq1} with $$A=\dg,\; B_1=\mb=\cg Du,\;B_2=d,\; G_1=g_2,\; G_2=cDu+f_2.$$ Then $ B_1^2,\mb^2\in \bbM(\Og,T)$ by a simple use of H\"older's inequality and the assumption that $\sup_{(0,T)}\|Du\|_{L^{2p}(\Og)}<\infty$. From the equation of $v$, since $u$ is already bounded, the conditions on $cDu,d,g_2,f_2$ and $\dg\in\bbM(\Og,T)$, we prove that $v$ is bounded and H\"older continuous under \mref{tri1} (this does not need differentiate the equation of $v$).

The condition \mref{triFa2} is exactly what we assumed in \mref{triFa} of \reftheo{genN}. The same argument, using \refrem{Duvrem}, there applies. So that, we obtain an estimate for $\sup_{(0,T)}\|Dv\|_{L^{2p}(\Og)}$.

We now obtain an estimate for $\sup_{(0,T)}\|DW\|_{L^{2p}(\Og)}$. Since $u,v$ are bounded so are the coefficients of the the system. The existence of $u,v$ on $(0,T)$ then follows from the result of \cite{Am2}. 
\eproof

\brem{Gap} The condition \mref{triFa2} was assumed because we differentiate the equation of $v$ as in the proof of \reftheo{genN}. Dropping \mref{triFa2},  we can still prove that $v$ is H\"older continuous, and so is $W$. If $\mA$ satisfies the spectral gap condition as in \cite{dleANS} the the argument there for systems also gives the existence of $W$ on $(0,T)$. \erem

\bibliographystyle{plain}

\end{document}